\documentclass[12pt]{article}
\usepackage{amsfonts,amssymb,amsmath}
\usepackage{graphicx}
\title{ Hamiltonian Cycle on the Set of Genotypes and Non-Ergodic Quadratic Stochastic Operators }
\author{Nasir N. Ganikhodjaev$^1,$ Uygun U. Jamilov$^{1,2}$ and Ramazon T. Mukhitdinov$^3$}
\begin{document}

\maketitle

\begin{center}
$^1$Department of Computational and Theoretical Sciences, Faculty
of Science, IIUM, 25200 Kuantan,Malaysia. nasirgani@hotmail.com\\
$^2$Institute of Mathematics,100125, Tashkent, Uzbekistan. jamilovu@yandex.ru\\
$^3$Bukhara Engineering-Technical Institute of High Technologies, 105017 Bukhara,Uzbekistan. ramazon-mukhitdinov@rambler.ru
 \vskip 0.1cm
\end{center}

\begin {abstract} On the set of genotypes $\Phi=\{1,\cdots,m\}$ we introduce  a binary relation  generated by Volterra quadratic stochastic operator $V$ on $(m-1)$ dimensional simplex $S^{m-1}$ and prove that
the operator $V$ be  non-ergodic if either there exists a Hamiltonian cycle or one of the vertices $M_i=(\delta_{1i},\delta_{2i},\cdots,\delta_{mi})$ of the simplex $S^{m-1}$ is a source and restriction of  $V$ to the invariant face $F_i=\{x\in S^{m-1}: x_i=0\}$ is non-ergodic. In this paper we prove this result for $m=2,3,4.$

\end{abstract}

{\bf Keywords}:  Volterra quadratic stochastic operator; binary relation; Hamiltonian cycle; non-ergodic transformation.

\section{Introduction}
Quadratic stochastic operators were first introduced by Bernstein \cite{Br}. Such operators frequently arises in many models of mathematical genetics, namely, theory of heredity \cite{Br},\cite{GRRU}-\cite{Lyu2}.
Consider a biological population,that is a community of
organisms closed with respect to reproduction. Assume that each
individual in this population belongs to precisely one species (genotypes)
 $1,\cdots, m . $ The scale of species is such that the species of the
parents $ i $ and $ j $ unambiguously determines the probability
of every species $k$ for the first generation of direct
descendants. Denote this probability, that is to be called the
heredity coefficient, by $p_{ij,k}.$ It is then obvious that
$p_{ij,k}\geq 0$ for all $i,j,k$ and that
$$\sum^m_{k=1}p_{ij,k}=1  (i,j,k=1,\cdots,m). $$

The state of the population can be described by the tuple
$(x_1,x_2,\cdots,x_m)$ of species probabilities,that is $ x_k$ is
the fraction of the species $k$ in the total population. In the
case of panmixia (random interbreeding) the parent pairs $i$ and
$j$ arise for a fixed state $\textbf{x}=(x_1,x_2,\cdots,x_m)$ with
probability $x_ix_j.$ Hence the total probability of the species
$k$ in the first generation of direct descendants is defined by
$$ \sum^m_{i,j=1}p_{ij,k}x_ix_j, \quad  (k=1,\cdots,m) $$

Let
\begin{equation} \label{simplex}
S^{m-1} =\{\textbf{x} = (x_1, x_2,\cdots,x_m)\in R^m : \mbox{ for any } i \ x_i \geq 0, \mbox { and } \sum_{i=1}^m x_i=1\}
\end{equation}
be the $(m-1)$-dimensional simplex.
A map $V$  of $S^{m-1}$  into itself is called a quadratic stochastic operator (q.s.o.) if
\begin{equation} \label{qso}
(V \textbf{x})_k = \sum_{i,j=1}^m p_{ij,k}x_ix_j
\end{equation}
 for any $\textbf{x} \in S^{m-1}$  and for all $ k = 1,\cdots,m,$  where
 $$a) p_{ij,k}\geq 0,
b) p_{ij,k} = p_{ji,k} \mbox { for all }  i, j, k; c) \sum_{k=1}^m p_{ij,k}=1.$$

Assume $\{\textbf{x}^{(n)}\in S^{m-1}: n=0,1,2,\cdots \}$ is the trajectory of the initial point $\textbf{x}\in S^{m-1},$
where $\textbf{x}^{(n+1)}=V(\textbf{x}^{(n)})$ for all $n=0,1,2,\cdots,$ with $ \textbf{x}^{(0)}=\textbf{x}.$

{\bf Definition 1.1 } A  point $\textbf{a}\in S^{m-1}$ is called a fixed point of a qso $V$ if $V(\textbf{a})=\textbf{a}.$

{\bf Definition 1.2} A qso $V$ is called regular if for any initial point $\textbf{x} \in S^{m-1}$  the limit
$$\lim_{n\rightarrow \infty}V^n(\textbf{x})$$
exists.

Note that the limit point be a fixed point of a qso $V.$ Thus the fixed points
of qso describe limit or long run behavior the trajectories of any initial point.
Limit behavior of trajectories and fixed points of qso play important role in many
applied problems \cite{GR1},\cite{GR2},\cite{J},\cite{K}, \cite{Lyu}.
The biological treatment of the regularity of qso is rather clear: in long run the
distribution of species in next generation coincide with distribution of species in
previous one, i.e., stable.

{\bf Definition 1.3 } A qso $V$ is said to be ergodic if the limit
\begin{equation}\label{erg}
\lim_{n\rightarrow \infty } \frac{1}{n} \sum_{k=0}^{n-1} V^k(x)
\end{equation}
exists for almost all  $\textbf{x}\in S^{m-1}$ with respect to usual Lebesgue measure on $S^{m-1}.$

On the basis of numerical calculations Ulam conjectured \cite{U} that
the ergodic theorem holds for any qso $V.$
  In 1977 Zakharevich \cite{Z} proved that this conjecture is false in general.
  Later in \cite{GZ} was established necessary and sufficient
condition to be non-ergodic transformation for qso defined on $S^2.$ The biological
treatment of non-ergodicity of qso $V$ is the following: in long run the behavior of distributions of species is chaotic, i.e., unpredictable.

Note that a regular qso $V$ is ergodic, but in generally from ergodicity does not follow regularity.
Let $\pi$ be a permutation of a set $\{1,2,\cdots, m\}.$ Then one can define
one-to-one transformation $T_\pi: S^{m-1}\rightarrow S^{m-1}$ as
\begin{equation}\label{per}
T_{\pi}(x_1,x_2,\cdots, x_m)=(x_{\pi(1)}, x_{\pi(2)}, \cdots,x_{\pi(m)})
\end{equation}
It is evident that a transformation $T_\pi^{-1}VT_\pi$ is qso and limit behavior
of trajectories of qso $T_\pi^{-1}VT_\pi$ coincide with limit behavior of trajectories of qso $V$ \cite{GE}. Thus for any permutation $\pi$ a qso $T_\pi^{-1}VT_\pi$ be regular (resp. ergodic) if and only if a qso $V$ is regular (resp. ergodic).

    The following notations will be used throughout this paper \cite{J}. We let $\partial S^{m-1}$ denote the boundary of $ S^{m-1}$:\\
  $$\partial S^{m-1}=\{\textbf{x}\in S^{m-1}: x_i=0 \mbox{ for at least one } i\in \{1,2,\cdots,m\}\},$$
  and also for $i=1,2,\cdots,m:$ \\
   $M_i=(\delta_{1i},\delta_{2i},\cdots,\delta_{mi})$ be $i$th vertex of the simplex $S^{m-1},$  \\
     $F_i$ be the $i$th face of $S^{m-1},$ where $F_i=\{\textbf{x}\in S^{m-1}: x_i=0\};$ \\
     and interior of $S^{m-1}$ be the set $intS^{m-1}=\{\textbf{x}\in S^{m-1}:x_1x_2\cdots x_m>0\}.$

     Also denote by  $\omega(x^0)$ $\omega$-set of limit points of the trajectory  $\{x^{(n)}, n=0,1,2,...\}.$

\section{Extremal Volterra Quadratic Stochastic Operators and Tournaments}
 Let $ \Phi=\{1,\cdots, m\} $ be the scale of genotypes and $V$ a quadratic stochastic operator .\\
  {\bf Definition 2.1 } The quadratic stochastic operator $V$ is called Volterra, if $p_{ij,k}=0 $ for any $k\notin \{i,j\}.$ \\
 The biological treatment of such operators is rather clear: the offspring repeats one of its parents. \\
 Evidently for any Volterra qso $p_{ii,i}=1$ for any $i=1,\cdots,m.$ A Volterra qso $V$ defined on $S^{m-1}$ has following form
 \begin{equation}\label{vol1}
(V\textbf{x})_k=x_k^2+2\sum_{i\neq k}^m p_{ik,k}x_ix_k,
\end{equation}
 where $k=1,\cdots,m.$

 {\bf Proposition 2.1 } \cite{GR1} \cite{J} A qso $V$ is a Volterra if and only if
 \begin{equation}\label{vol}
(V\textbf{x})_k=x_k(1+\sum_{i=1}^ma_{ki}x_i)
\end{equation}
where $A=(a_{ij})_1^m$ is a skew-symmetric matrix with $a_{ki}=2p_{ik,k}-1,$ and $|a_{ij}|\leq 1.$ Here $i,j\in \{1,2,\cdots,m\}.$\\
{\bf Proposition 2.2 } Let $V$ be a Volterra qso. Then for any permutation $\pi \in S_m$ the transformation $T_\pi^{-1} V T_\pi $
is the Volterra qso. \\
{\it Proof } For qso $V$ (5) and permutation $\pi \in S_m$ simple algebra gives
 \begin{equation}\label{vol11}
(T_\pi^{-1} VT_\pi\textbf{x})_k=x_k^2+2\sum_{i\neq k}^m p_{\pi^{-1}(i)\pi^{-1}(k),\pi^{-1}(k)}x_ix_k,
\end{equation}
 where $k=1,\cdots,m.$ It is evident that
the transformation (7) is the Volterra qso.

 In \cite{GR1} and \cite{GR2} the theory of Volterra qso was developed using theory of the Lyapunov functions and tournaments. Below we will consider a class of so-called extremal qso.\\
 {\bf Definition 2.2 } The quadratic stochastic operator $V$ is called extremal Volterra, if $p_{ij,k}=0  \mbox{ or } 1 $ for  $k \in \{i,j\}.$ \\
Let $\mathcal{E}$ be the set of all extremal Volterra qso on $S^{m-1}.$ Since $p_{ik,i}=1-p_{ik,k},$ the total number of mixed products $x_ix_j$ with $i\neq j$ in (5) is equal to $m(m-1)/2.$ Note that if $V$ be an extremal Volterra qso then elements of matrix $A$ take only two values $\pm 1.$ Thus there are $2^{\frac{m(m-1)}{2}}$ extremal Volterra qso. \\
{\bf Proposition 2.3 } Let $V$ be an extremal Volterra qso. Then for any permutation $\pi \in S_m$ the transformation $T_\pi^{-1} V T_\pi $
is the extremal Volterra qso. \\
{\it Proof } The proof immediately follows from Proposition 2.2.  \\
{\bf Definition 2.3 } Let $V_1,V_2\in \mathcal{E}$ be  two extremal Volterra qso.  We will say that the qso $V_1$ is equivalent to $V_2$ and denote $V_1\sim V_2$ if there exists a permutation $\pi\in S_m$ such that $V_2=T_\pi^{-1} V_1 T_\pi .$ \\
Since this relation is equivalence one we can part the set $\mathcal{E}$ into equivalence classes.
Let  $V$ be extremal Volterra quadratic stochastic operator. Then on the set  $ \Phi=\{1,\cdots, m\} $ of genotypes  one can introduce following a binary relation: if $ p_{ij,i}=1$  we will say that the genotype $i$ \emph{dominate} genotype $j$  and denote $i\succ j.$
 Any two genotypes $i$ and $j$ are comparable,i.e.,  $i\succ j $ or $j\succ i. $  Let the set $ \Phi=\{1,\cdots, m\} $ is provided with directed graph structure where the edge connecting genotypes $i$ and $j$ directed from $i$ to $j$ if $i\succ j.$ That is, it is a directed complete graph in which every pair of vertices is connected by a single directed edge. Such graphs is called a tournament. \\
{\bf Definition 2.4 } A  cycle of a graph $ \Phi$ is a subset of the edge set of $ \Phi$ that forms a path such that the first node of the path corresponds to the last. \\
{\bf Definition 2.5 } A cycle that uses graph vertex of a graph exactly once is called a Hamiltonian cycle. \\
 A graph containing no cycles of any length is known as an acyclic graph, whereas a graph containing at least one cycle is called a cyclic graph.\\
 Let the  indegree of a vertex is the number of edges leading to that vertex, and the outdegree of a vertex is the number of edges leading away from that vertex. A vertex with an indegree of $0$ is called a source (since one can only leave it) and a vertex with an outdegree of $0$  is called a sink (since one cannot leave it). \\
 {\bf Proposition 2.4 } If for given qso $V$ a vertex $M_i=(\delta_{1i},\delta_{2i},\cdots,\delta_{mi})$ of the simplex $S^{m-1}$ is a source or sink, then the face $F_i=\{x\in S^{m-1}: x_i=0\}$ is invariant subset. \\
{\it Proof } The proof immediately follows from definition of source and sink.  \\
 \section{Ergodicity of Extremal Volterra Operators }

 Let $\Phi= \{1,\cdots,m\}$ be a set of genotypes. In this paper we consider small $m$ with $m=2,3,4.$
 General case will be subject of next paper.\\
 {\bf Theorem 3.1 } An extremal Volterra quadratic stochastic operator $V$ be non-ergodic if either there exists a Hamiltonian cycle or one of the vertices $M_i=(\delta_{1i},\delta_{2i},\cdots,\delta_{mi})$ of the simplex $S^{m-1}$ is a source and restriction of  $V$ to the invariant face $F_i=\{x\in S^{m-1}: x_i=0\}$ is non-ergodic.

 {\it Proof }
 For $k=2$ we have two extremal Volterra qso
 \begin{equation} \label{operator4}
 \begin{array}{llll}
x^\prime_1= x^2_1 + 2\alpha x_1x_2,\\[2mm]
x^\prime_2= x^2_2+2(1-\alpha) x_1x_2 \\[2mm]
\end{array}
\end{equation}
where $\alpha \in \{0,1\}$. Since we have two genotypes only, the corresponding graph is acyclic
and simple analysis shows that the operator (8) is the regular. It is evident that the qso (8) has
two fixed points $M_1=(1,0)$ and $M_2=(0,1).$ For $\alpha=0$ (respectively $\alpha=1$) any trajectory converges to $M_2$ (respectively $M_1$).

Let $k=3.$ Then we have $8$ extremal Volterra operators
\begin{equation} \label{operator5}
 \begin{array}{llll}
x^\prime_1= x^2_1 + 2\alpha x_1x_2+2\beta x_1x_3,\\[2mm]
x^\prime_2= x^2_2+2(1-\alpha) x_1x_2+2\gamma x_2x_3 \\[2mm]
x^\prime_3= x^2_3+2(1-\beta) x_1x_3+2(1-\gamma) x_2x_3 \\[2mm]
\end{array}
\end{equation}
where $\alpha,\beta,\gamma \in \{0,1\}$.
It is easy to see that there exists Hamiltonian cycle with respective binary relation  generated by
qso with  $(\alpha,\beta,\gamma)=(0,1,0)$ or   $(\alpha,\beta,\gamma)=(1,0,1)$ and for remaining 6 cases
corresponding graph is acyclic. The qso (9) with $(\alpha,\beta,\gamma)=(1,0,1)$ has following form:
\begin{equation} \label{operator6}
 \begin{array}{llll}
x^\prime_1= x^2_1 + 2x_1x_2,\\[2mm]
x^\prime_2= x^2_2+2x_2x_3 \\[2mm]
x^\prime_3= x^2_3+2x_1x_3 \\[2mm]
\end{array}
\end{equation}
This operator was considered by Zakharevich \cite{Z} and was proven that it is non-ergodic transformation.
Second operator with $(\alpha,\beta,\gamma)=(0,1,0)$ is reduced to (10) by a permutation of genotypes.
For $(\alpha,\beta,\gamma)=(0,0,0)$ the corresponding qso has following form:
\begin{equation} \label{operator7}
 \begin{array}{llll}
x^\prime_1= x^2_1 ,\\[2mm]
x^\prime_2= x^2_2+2x_1x_2 \\[2mm]
x^\prime_3= x^2_3+2x_1x_3+2x_2x_3. \\[2mm]
\end{array}
\end{equation}
and all others qso with acyclic graph are reduced to (11) by some permutation of genotypes. Note that
there exist exactly 6 permutations and respectively 6 qso with acyclic graph. Simple analysis shows that qso (11)
is recurrent.One can rewrite qso (11) as following:
\begin{equation} \label{operator17}
 \begin{array}{llll}
x^\prime_1= x^2_1 ,\\[2mm]
x^\prime_2= x^2_2+2x_1x_2 \\[2mm]
x^\prime_3= x_3(2-x_3). \\[2mm]
\end{array}
\end{equation}
It is evident that one-dimensional transformation $\varphi(x)=x(2-x)$ has two fixed points $x^*=0$ and $x^{**}=1$ with $\varphi^\prime(0)>1$
and $\varphi^\prime(1)<1,$ such that $x^*$ is repelling (resp. $x^{**}$ attracting ) fixed point.
Since $\varphi(x)$ is increasing function on segment $[0,1]$ and $x^{**}$ is fixed point, then
$\lim_{n\rightarrow \infty}\varphi^{(n)}(x) =1,$ where $\varphi^{(n)}(x)=\varphi(\varphi^{(n-1)}(x)),$
therefore any trajectory of qso (12) converges to
vertex $(0,0,1),$ except fixed point $(1,0,0).$

  For $k=4$ the set $\mathcal{E}$ consists of 64 extremal Volterra operators
 \begin{equation} \label{operator8}
 \begin{array}{llll}
x^\prime_1= x^2_1 + 2\alpha x_1x_2+2\beta x_1x_3+2\gamma x_1x_4,\\[2mm]
x^\prime_2= x^2_2+2(1-\alpha) x_1x_2+2\delta x_2x_3+2\varepsilon x_2x_4 \\[2mm]
x^\prime_3= x^2_3+2(1-\beta) x_1x_3+2(1-\delta) x_2x_3+2\lambda x_3x_4 \\[2mm]
x^\prime_4= x^2_4+2(1-\gamma) x_1x_4+2(1-\varepsilon) x_2x_4+2(1-\lambda) x_3x_4 \\[2mm]
\end{array}
\end{equation}
where $\alpha,\beta,\gamma,\delta,\varepsilon, \lambda \in \{0,1\}$.
Simple but tedious algebra gives that the set $\mathcal{E}$ is parted into $4$ equivalent classes.
First class consists of $24$ extremal Volterra operators such that there exists Hamiltonian cycle.  Second class consists of $8$  extremal Volterra operators such that for corresponding graph one vertex is a sink and other $3$ vertices forms three cycle. Third class consists of $8$  extremal Volterra operators such that for corresponding graph one vertex is a source and other $3$ vertices forms three cycle. Finally, fourth class consists of $24$ extremal Volterra operators such that the corresponding graphs are acyclic.
Note that for $k=4$ there exist exactly 24 permutations.

\subsection{ First Equivalent Class}
If $\alpha=\beta=\delta=\varepsilon=\lambda=1$ and $\gamma=0$ then there exists the Hamiltonian
cycle $1\succ 2\succ3 \succ 4 \succ 1$ and corresponding qso (13) has following form
\begin{equation} \label{operator9}
 \begin{array}{llll}
x^\prime_1= x^2_1 + 2x_1x_2+2x_1x_3,\\[2mm]
x^\prime_2= x^2_2+2x_2x_3+2x_2x_4 \\[2mm]
x^\prime_3= x^2_3+2x_3x_4 \\[2mm]
x^\prime_4= x^2_4+2x_1x_4. \\[2mm]
\end{array}
\end{equation}
The equivalent class containing this operator consist of $24$ extremal Volterra qso and for each of them corresponding tournament contains the Hamiltonian cycle.
It is easy to see that the set of fixed points contains following points $M_1(1,0,0,0),$
$M_2(0,1,0,0),$ $M_3(0,0,1,0), M_4(0,0,0,1), M_5(1/3,0,1/3,1/3),M_6(1/3,1/3,0,1/3).$ Introducing new variables $y_1=x_1, y_2=x_2+x_3, y_4=x_4,$ one can see that
 $$y^\prime_1=x^\prime_1, y^\prime_2=x^\prime_2+x^\prime_3, y^\prime_4=x^\prime_4,$$ and qso (14) reduced to following
 \begin{equation} \label{operator91}
 \begin{array}{llll}
y^\prime_1= y^2_1 + 2y_1y_2,\\[2mm]
y^\prime_2= y^2_2+2y_2y_3 \\[2mm]
y^\prime_3= y^2_3+2y_1y_3 \\[2mm]
\end{array}
\end{equation}
Note that this qso (\ref{operator9}) has following invariant subset $\{\tau M_5+(1-\tau)M_6, \ \ \tau \in [0,1] \}$ and trajectory of (\ref{operator9})  with initial point from
these invariant subset converges to fixed point $M_6$. Since Lebesgue measure of invariant subset is equal to $0$, this qso is non-ergodic. Thus we have proved following Proposition. \\
{\bf Proposition 3.1 } Any qso from first class is non-ergodic operator.

\subsection{ Second Equivalent Class}
If $\alpha=\beta=\gamma=\delta=\lambda=1$ and $\varepsilon=0$ then the vertex $M_1=(1,0,0,0)$ is the sink, other three vertices forms  three cycle $2\succ 3\succ 4\succ 2$ and corresponding qso (13) has following form
\begin{equation} \label{operator18}
 \begin{array}{llll}
x^\prime_1= x^2_1 + 2x_1x_2+2x_1x_3+2 x_1x_4,\\[2mm]
x^\prime_2= x^2_2+2x_2x_3 \\[2mm]
x^\prime_3= x^2_3+2x_3x_4 \\[2mm]
x^\prime_4= x^2_4+2x_2x_4. \\[2mm]
\end{array}
\end{equation}
The equivalent class containing this operator consist of $8$ extremal Volterra qso and for each of them corresponding tournament contains three cycle. It is easy to see that the set of fixed points consists of following points $M_1(1,0,0,0),$ $M_2(0,1,0,0),$ $M_3(0,0,1,0),$ $M_4(0,0,0,1), M_5(0,1/3,1/3,1/3).$ First equality in (16) one can rewrite as $x_1^\prime=x_1(2-x_1)$ similar (12). As proved above, for any initial $x^0\in \partial S^3$ its trajectory converges to vertex $M_1=(1,0,0,0).$ It is evident that the face $F_1=\{x\in S^3:x_1=0\}$
is invariant subset with $0$ Lebesgue measure and restriction of qso (16) to $F_1$ is reduced to Zakharevich's example (10).Thus we have proved following Proposition

{\bf Proposition 3.2 } Extremal Volterra  qso $V$ belonging  to the second class is ergodic.

\subsection{ Third Equivalent Class}
If $\alpha=\beta=\delta=\varepsilon=1$ and $\gamma=\lambda=0$ then the vertex $M_3=(0,0,1,0)$ is a source , other three vertices forms three cycle $1\succ 2\succ 4\succ 1$ and corresponding qso (13) has following form
\begin{equation} \label{operator181}
 \begin{array}{llll}
x^\prime_1= x^2_1 + 2x_1x_2+2x_1x_3,\\[2mm]
x^\prime_2= x^2_2+2x_2x_3+2x_2x_4 \\[2mm]
x^\prime_3= x^2_3 \\[2mm]
x^\prime_4= x^2_4+2 x_1x_4+2x_3x_4\\[2mm]
\end{array}
\end{equation}
The equivalent class containing this operator consist of also $8$ extremal Volterra qso and for each of them corresponding tournament contains three cycle. It is easy to see that the set of fixed points consists of following points $M_1(1,0,0,0),$ $M_2(0,1,0,0),$ $M_3(0,0,1,0),$ $M_4(0,0,0,1), M_5(1/3,1/3,0,1/3).$ It is evident that  $\lim\limits_{ n \rightarrow \infty}x_3^{(n)}=0$. If for initial $x^0\in S^3$ we have $x_1^0\neq 0$ then  $\omega(x^0) \subset F_3=\{x\in S^3:x_3=0\}$.  It is evident that the face $F_3=\{x\in S^{m-1}:x_3=0\}$
is invariant subset and restriction of qso (17) on $F_3$ is reduced to Zakharevich's example (10).Thus we have proved following Proposition

{\bf Proposition 3.3 } If extremal Volterra  qso $V$ belongs to the third class, then it is non-ergodic transformation.
\subsection{ Fourth Equivalent Class}
If $\alpha=\beta=\gamma=\delta=\varepsilon=\lambda=1$  then corresponding graph is acyclic and  corresponding qso (13) has following form
\begin{equation} \label{operator28}
 \begin{array}{llll}
x^\prime_1= x^2_1 + 2 x_1x_2+2x_1x_3+2 x_1x_4,\\[2mm]
x^\prime_2= x^2_2+2 x_2x_3+2 x_2x_4 \\[2mm]
x^\prime_3= x^2_3+2x_3x_4 \\[2mm]
x^\prime_4= x^2_4 \\[2mm]
\end{array}
\end{equation}
The equivalent class containing this operator consist of  $24$ extremal Volterra qso and for each of them corresponding tournament is acyclic. It is easy to see that the set of fixed points consists of following points $M_1(1,0,0,0),$ $M_2(0,1,0,0),$ $M_3(0,0,1,0), M_4(0,0,0,1)$
and the sets $F_1$ and $\Gamma_{34}=\{\textbf{x}\in S^3: x_1=x_2=0\}$ are invariant. First equality  (18) one can rewrite as $x_1^\prime=x_1(2-x_1)$ and using properties of the function $\varphi(x)=x(2-x)$
 we have \\
1) If $x_1^0\neq 0$ then any trajectory of qso (18) converges to the vertex $M_1=(1,0,0,0);$ \\
2) if an initial point belongs to set  $F_1\setminus \Gamma_{34},$ i.e.,  $x_1^0= 0,$ then the operator (18) is reduced to qso (11),
 therefore any trajectory converges to the vertex $M_2=(0,1,0,0);$ \\
 3) if $\textbf{x}\in \Gamma_{34}\setminus M_4 $ then the trajectory converges to third vertex $M_3=(0,0,1,0).$
 Thus we have the following Proposition.

{\bf Proposition 3.4 } Any extremal Volterra  qso $V$ belonging  to the fourth class is regular transformation.

 Collecting together all four Propositions we complete the proof of Theorem 3.1

 \section{ Conclusion }

 Above we consider extremal Volterra qso on the set $\Omega$ with $|\Omega|=m$ where we limit ourself by $m=2,3,4.$  To prove this result for arbitrary $m$ we need to apply more deep properties of the tournament and it cycles.

\vskip 0.6 truecm

 {\bf Acknowledgment.} The second-named author, Post Doctorate Researcher (U.U.J) thanks International Islamic University Malaysia for kind hospitality and providing all facilities.

 \vskip 0.3 truecm

\begin{center}

\end{center}
\end{document}